\documentclass[12pt]{article}
\usepackage{amsfonts}
\normalbaselines
\catcode `\@=11
\@addtoreset{equation}{section}

\def\section{\@startsection {section}{1}{\z@}{-3.5ex plus -1ex minus
     -.2ex}{2.3ex plus .2ex}{\normalsize\bf}}
\def
\subsection{\@startsection{subsection}{2}{\z@}{-3.25ex plus -1ex minus
 -.2ex}{1.5ex plus .2ex}{\normalsize\bf}}
\def\thebibliography#1{\section*{References\markboth
  {REFERENCES}{REFERENCES}}\list
  {[\arabic{enumi}]}{\settowidth\labelwidth{[#1]}\leftmargin\labelwidth
  \advance\leftmargin\labelsep
  \usecounter{enumi}}
  \def\newblock{\hskip .11em plus .33em minus -.07em}
  \sloppy
  \sfcode`\.=1000\relax}
 
\catcode `\@=12
\newcommand{\gk}{\mbox{\ai K}}%K mathematical italic cmbsy for Hilbert space
\newcommand{\df}{\mbox{$:=$}}% def
\newtheorem{rem}{Remark}
\newtheorem{com}{Comment}

\newtheorem{thm}{Theorem}
\newtheorem{pr}{Proposition}
\newtheorem{cor}{Corollary}

\font\aa=msam10 scaled\magstep1%mat,
\font\aba=msbm7 % for C, Z, R in small format
\font\ai=eusb10 scaled\magstep1% script bold
\newcommand{\prend}{\aa \symbol{'003}}%symbol for end of the proof
\newcommand{\pre}{\mbox{\prend}}
\newcommand{\gh}{\mbox{\ai H}}%H mathematical italiac cmbsy for Hilbert space
\newcommand{\gph}{\mbox{{\bf P}(\gh )}}% Projective Hilbert space PH
\newcommand{\gz}{\mbox{$\mathbb Z$}}%Z double
\newcommand{\gcp}{\mbox{${\mathbb{CP}}$}}%CP double
\newcommand{\gc}{\mbox{$\mathbb C$}}
\newcommand{\gcm}{\mbox{\aba C}}% C double for exponent
\newcommand{\Ka}{K\"ahler}
\newcommand{\Gras}{$\mbox{$G_n({\mathbb C}^{m+n})$}$}% Grassmannian
\newcommand{\men}{\mbox{$\widetilde {\bf M}$}} %Manifold
\newcommand{\got}[1]{{\mbox{${\mathfrak{#1}}$}}}% gothic with mbox for math.
\newcommand{\mb}[1]{{\mbox{\boldmath{$#1$}}}}% mathematical bold

\begin{document}

\vspace*{2.5cm}
\noindent
{ \bf  COHERENT STATES, LINE BUNDLES AND DIVISORS}\vspace{1.3cm}\\
\noindent
\hspace*{1in}
\begin{minipage}{13cm}
Stefan Berceanu$^{1}$ \vspace{0.3cm}\\
 $^{1}$ Institute for Physics and Nuclear Engineering\\
 \makebox[3mm]{ }Department of Theoretical Physics\\
 \makebox[3mm]{ }PO BOX MG-6, Bucharest-Magurele,\\
 \makebox[3mm]{ }Romania; E-mail: Berceanu@theor1.ifa.ro \\
\end{minipage}

\vspace*{0.5cm}
  \begin{abstract}
\noindent
For homogeneous simply connected  Hodge manifolds it is proved that the set
of coherent vectors orthogonal to a given one is the divisor responsible for
the homogeneous holomorphic line bundle of the coherent vectors. In particular,
for naturally reductive spaces, the divisor is the cut locus.
\end{abstract}

\section{\hspace{-4mm}.\hspace{2mm}INTRODUCTION}

The coherent states \cite{klauder,per,raw}
are a powerful tool in global differential geometry\cite{ber1,ber3}.
For example, the remark  {\bf Polar divisor = Cut locus}, proved
on naturally reductive spaces \cite{ber2}, gives a description of the
cut locus in terms of coherent states.	Let us briefly recall these
notions.

 Let $X$ be complete Riemannian manifold. The point  $q$ is in the
{\it  cut locus} ${\bf CL_p}$ of $p\in X$ if $q$ is the nearest point to $p$ on
the geodesic emanating from $p$ beyond which the geodesic ceases to minimize
his arc length (cf. \cite{kn}, see also Ref. \cite{ber2} for more references).

We call {\it  polar divisor} of $\mb{e}_0$ the set
$\Sigma _0=\left\{ \mb{e} \in e(G)|(\mb{e}_0, \mb{e})=0\right\}$, where
$e(G)$ is the family of coherent vectors \cite{raw}. This
denomination is inspired after Wu \cite{wu}, who used this term in the case
of the complex Grassmann manifold.

In this paper we shall emphasize an aspect of the  deep relationship
between coherent states and algebraic geometry.  Indeed, the notion of polar
divisor, introduced in the context of coherent states, is in agreement with the
notion of divisor in algebraic geometry \cite{hartshorne}.
The main result of
this paper is the establishment of the relationship between vector
coherent state manifold ${\bf M}$, viewed as a holomorphic homogeneous
line bundle over the coherent state manifold \men , and the polar divisor. The
set $\Sigma_0$ is the divisor responsible for the line bundle ${\bf M}$.  The
result is proved for homogeneous simply connected Hodge manifolds. In
particular, for naturally reductive spaces, the divisor responsible for
${\bf M}$ is the cut locus ${\bf CL_0}$. The enunciation of this theorem was
included in Ref.
\cite{sberezin}. Using Rawnsley's definition of coherent states
\cite{raw} instead of Perelomov's one, it is possible \cite{sbs}
to drop out the restriction that the manifold $\men$ to be a homogeneous one.

\newpage
The lay out of this paper is as follows: \S 2 collects some feature needed
during
the paper -- a short remember of the notion of coherent states viewed as
homogeneous line bundles, a breviary on divisors, and
a brief review of the results established in Ref. \cite{ber2} on cut locus
and coherent states. The main results are proved in \S 3. The illustration on
the complex Grassmann manifold follows the notation of \cite{ber4}.

\section{\hspace{-4mm}.\hspace{2mm}INGREDIENTS}
\subsection{\hspace{-5mm}.\hspace{2mm}Homogeneous Line Bundles and
Coherent States}

1.)  Let us consider the principal bundle
\begin{equation}
\label{pbundle}
K\stackrel{i}{\rightarrow} G\stackrel{\lambda}{\rightarrow}\men,
\end{equation}
where \men~ is diffeomorphic with $G/K$, $i$ is the inclusion and
$\lambda $ is the natural projection $\lambda (g)=gK$. Let $\chi$ be
 a continuous representation of the group $K$ on the Hilbert space
\gk~ and 
 let
${\bf M}_{\chi}\df
\men \times_{\chi}\gk$, or simply ${\bf M}\df\men \times_K\gk$, be
the $G$-{\it homogeneous vector bundle\/} \cite{bott}
associated by
 $\chi$  to the principal $K$-bundle (\ref{pbundle}).
Let  $U\subset\men$~ be open. We introduce the notation
\begin{equation}
\label{local}
(G)^{U}=\{g\in G| go\in U\},
\end{equation}
where $o$ is  the base point in  $ \men$.
Then the continuous (holomorphic) sections of
${\bf M}_{\chi}$ over $U$ are precisely the continuous (resp. holomorphic)
maps $\sigma :U\rightarrow G\times_{\chi}\gk$ of the form
\begin{equation}
\label{sectiuni}
\sigma(go)=[g,e_{\sigma}(g)],~e_{\sigma}:(G)^U
\rightarrow \gk,
\end{equation}
where $e_{\sigma}$ satisfies the ``functional equation'':
\begin{equation}
\label{fec}
e_{\sigma}(gp)=\chi (p)^{-1}e_{\sigma}(g), g\in(G)^U, p\in K.
\end{equation}

2.) Let $ \xi :\gh^{\star}=\gh \setminus \{0\}\rightarrow\gph,
~ \xi (\mb{z})=[\mb{z}]$
be the mapping	which associates to the point $\mb{z}$ in the punctured
Hilbert space the linear
subspace $[\mb{z}]$ generated by $\mb{z}$, where $[\lambda\mb{z}]=[\mb{z}],
\lambda\in\gc^{*}$.% The hermitian scalar product $(\cdot,\cdot )$ on \gh~
%is linear in the second argument.

 Let us consider the principal bundle (\ref{pbundle}) and let us suppose
the existence of a map $e:G\rightarrow \gh^*$  as in eq. (\ref{sectiuni})
with the property (\ref{fec}) but globally
defined, i.e. on the neighbourhood (\ref{local}) $(G)^{\men}$.
Then $e(G)$ is called {\it  family of coherent vectors\/} \cite{raw}.
If there is a morphism of principal bundles,
i.e. the following diagram is commutative,
\begin{eqnarray}
\label{coherent}
\nonumber G &
~\stackrel{e}{\longrightarrow}~
& \gh^\star \\
\lambda\downarrow &  & \downarrow\xi\\ \label{comdiag}
\nonumber\men & ~\stackrel{\iota}{\longrightarrow}~ & \gph
\end{eqnarray}
then $\iota (\men )$  is called {\it family of coherent states corresponding
to the family of coherent vectors} $e(G)$ \cite{raw}.
 The manifold \men ~ is called
{\it coherent state manifold} and the
$G$-homogeneous line bundle ${\bf M}_{\chi}$ is called {\it coherent vector
manifold} \cite{sbcag}.

We impose the following restrictions:
a) the mapping $\iota$	is an embedding
in some projective Hilbert space
\begin{equation}\label{emb}
\iota :
\widetilde {\bf M} \hookrightarrow \gph ;
\end{equation}
b) \men ~ is homogeneous \cite{per};
c)  the embedding $\iota$ is  k\"ahlerian (cf. Ch. 8 in
\cite{besse}); d) the  line bundle ${\bf M}$ is {\it very ample}. If
$\chi$ induces a unitary representation $\pi_{\chi}$ on the Hilbert
space of holomorphic sections $\gh = \Gamma_{{\mbox{{\rm {hol}}}}}(\men ,{\bf
M}_{\chi}) =H^0(\men , {\bf M}_{\chi} )$, then we have an effective
realization of Perelomov's coherent states.

The  Perelomov's coherent vectors  are
\begin{equation}
{\mb e}_{Z,j}=(\exp\sum_{{\varphi}\in\Delta^+_n}(Z_{\varphi}F^+_{\varphi}))
j ,~~~~\underline{{\mb e}}_{Z,j}=({\mb e}_{Z,j} ,
{\mb e}_{Z,j})^{-1/2}{\mb e}_{Z,j}.
\label{z}
\end{equation}
In eq. ({\ref{z})
$\Delta^+_n$ denotes the positive non-compact roots, $ Z\df (Z_ \varphi )
\in {\gc}^n$  are local
coordinates in the  neighbourhood
${\cal V}_0 \subset \widetilde {\bf M} $ of the base point $o$,
$F^+_{\varphi} j\not= 0,  ~
(F^-_{\varphi} j = 0),$
$ \varphi \in\Delta^+_n$,  and $j$ is the dominant
weight vector of the representation. The system $\{e(g)\}$,
$g\in G^{\gcm}$ is overcomplete \cite{per,onofri} and $(e(g),e(g'))$, up to a
factor, is a reproducing kernel for the holomorphic line bundle
 ${\bf M_{\chi}}\rightarrow\men$ \cite{zb}.

\subsection{\hspace{-5mm}.\hspace{2mm}Divisors}

Let \men~ be a	complex manifold and $D$ a
locally finite formal combination $D=\sum a_iV_i$ of irreducible
analytic not necessarilly smooth  hypersurfaces of \men .  $D$ is called
a {\it Weil
divisor} of \men .
Let ${\mbox{Div}}(\men )$ denote the abelian group  of Weil divisors. If
$a_i\ge 0$,  then $D$ is  an
{\it effective} ({\it holomorphic})  divisor and this is denoted $D\ge 0$.

The divisors can be defined   in terms of sheaf theory.
A {\it Cartier divisor} $D$  on an algebraic variety \men~ is a
global section of  the quotient sheaf ${\cal D}$. We use the following
notation:
${\cal A}$
-- the sheaf of germs of complex $C^{\infty}$-valued functions;
${\cal A}^*$ -- the sheaf of germs of complex $C^{\infty}$-valued functions
nowhere zero;
${\cal O}$ -- the sheaf of germs of holomorphic functions;
${\cal O}^*$  --  the sheaf of germs holomorphic functions vanishing nowhere;
${\cal M}$ -- the sheaf of germs of local meromorphic  functions;
${\cal M}^*$ -- the sheaf of germs of invertible meromorphic functions;
${\cal D}={\cal M}^*/{\cal O}^*$.

Note that {\it the equivalence classes of $C^{\infty}$ line
bundles are in 1-1 correspondence with the elements of the cohomology group
$H^1(\men ,{\cal A}^*)$} and there is an
{\it isomorphism between the group of continuous line bundles and the second
cohomology group with integer coefficients}:
$ H^1(\men ,{\cal A}^*)\stackrel{\delta}{\rightarrow}H^2(\men ,\gz )$.
However, for holomorphic line bundles,
there is  only a homomorphism
$H^1(\men ,{\cal O}^*)\stackrel{{\delta}^1}{\rightarrow}H^2(\men ,\gz )$.

On complex manifolds there is an {\it the isomorphism of Weil and Cartier
groups of divisors}:
${\mbox{Div}}(\men )=H^0(\men ,{\cal D}).$

There is a  functorial homomorphism $[~]$ between the group of divisors and
the Picard group of  equivalence class of  holomorphic line bundles
$[~]:{\mbox{Div}}(\men ){\rightarrow} H^1(\men ,{\cal O}^*)$.
$[D]$ is a $\gc^*$-bundle, but we denote by the same symbol an analytic line
bundle determined up to an isomorphism.

The exact sequence
$1\rightarrow{\cal O}^*\stackrel{i}{\rightarrow}
{\cal M}^*\stackrel{k}{\rightarrow}{\cal D}\rightarrow 1$
induces the exact cohomology sequence
\begin{equation}\label{exseq}
0\rightarrow H^0(\men ,{\cal O}^*)\stackrel{i^0}{\rightarrow}
H^0(\men ,{\cal M}^*)
\stackrel{k^0}{\rightarrow}H^0(\men ,{\cal
D})\stackrel{\delta^0}{\rightarrow}H^1(\men ,{\cal O}^*)\rightarrow\cdots
\end{equation}
The quotient group
${\mbox {\rm Cl}}(\men )\df H^0(\men ,{\cal D})/k^0H^0(\men ,{\cal M}^*)$,
called the {\it  group of divisor cla\-sses with respect to linear equivalence},
is isomorphic to a subgroup of $H^1(\men ,{\cal O}^*)$.
If  $f$ is a meromorphic function on \men , then
  $(f)$ is its associated divisor.
The divisors $D,~D'$ are {\it linearly equivalent} if $D=D'+(f)$,
where $f\in{\cal M}^*(\men )$.

The holomorphic {\it divisor} $D$ is said  {\it non-singular}
if, with respect to some open covering
$ U=\{U_i\}$, it is represented by place functions $f_i$ with the
property:
either $f_i\equiv 1$ or $U_i$ admits a system of local complex coordinates for
which $f_i$ is one of the coordinates.

\newpage
\subsection{\hspace{-5mm}.\hspace{2mm} Cut Locus}

\begin{rem} $ {\mbox{\rm codim}}_{\gcm} {\bf CL}_p \geq 1$.
\label{rem4}
\end{rem}
Let ${\got g}$ be the Lie algebra of the group $G$.
Let  ${\got g}={\got k}\oplus {\got m}$ be  the orthogonal
decomposition with respect to the $B$-form as explained below,
${\rm Exp}_p: T_p\widetilde {\bf M}\rightarrow \widetilde {\bf M}$  the
geodesic exponential map from the tangent space to the manifold  and
 $\exp :{\got g} \rightarrow G $ the exponential from the Lie algebra
to the group.

 Let us consider  the conditions:
\begin{itemize}
\item[A)]
 ${\rm Exp}\vert _o=\lambda \circ \exp \vert _{\got m} .$
\item[B)] On the Lie algebra {\got g} of $ G$ there exists and
$ Ad(G)$-invariant, symmetric, non-degenerate bilinear form $ B$ such that
the restriction of $B$ to the Lie algebra {\got k}  of $ K$  is likewise
non-degenerate.
\end{itemize}
Note that {\it for  the homogeneous space} $\widetilde {\bf M}\approx
G/K$ {\it properety } $B)$, {\it implies} $A)$.
The  condition $A)$ is verified by the symmetric spaces, but also   by the
naturally reductive spaces because they verify the condition  $B)$ .

\begin{thm}\label{cutu}
Let \men~   be a homogeneous manifold
$\widetilde {\bf M}\approx G/K$.  Suppose that there exists a unitary
irreducible representation $\pi_j$ of $G$ such that in a neighbourhood
$ {\cal V}_0$  around $Z=0$ the coherent states are parametrized as
in eq. (\ref{z}).  Then the manifold \men~ can be represented as
the disjoint union
$\widetilde {\bf M} ={\cal V}_0\cup \Sigma _0.$
Moreover, if the condition  $B)$ is true, then
$\Sigma _0={\bf CL}_0. $
\end{thm}
%%%%

\begin{cor}\label{cul}
Suppose that $\widetilde {\bf M}$ verifies  $ B)$  and admits the embedding
(\ref{emb}). Let $0, Z\in
\men $.  Then  $Z\in {\bf CL}_0$ iff
$d_c(\iota (0),\iota (Z))=\pi /2$, where
$d_c([{\mb{\omega}} '],[{\mb{\omega}} ])
=\arccos \frac{\vert({\mb{\omega}} ',{\mb{\omega}})\vert}
{\Vert {\mb{\omega}} '\Vert \Vert {\mb{\omega}} \Vert }~$. \end{cor}

We remember the explicit expression
of the cut locus on the complex projective space and Grassmannian.
\begin{rem}
{\it  On $\gcp^n$,  ${\bf CL}_0=\Sigma_0=H_1=\gcp^{n-1}$}.
\end{rem}
{\it Proof.}
Let the notation
 $ {\cal V}_i=\{z\in \gh^*|z_i\not= 0\}$, ${\cal U}_i=\xi ({\cal V}_i)$,
$H_i=\gph\setminus
{\cal U}_i$. The point
$p_0=[(1,0,0,\ldots  0)]\in \gph$ corresponds to the point
$0$ in the Remark. Then
the solution of the equation $([p_0],[z])=0$ is $[z]=[(0,\times ,\times
,\ldots  \times)]= H_1=\gcp^{n-1}\subset\gcp^{n} $ for $\gh ={\gc }^{n+1}$.\pre

The complex Grassmann manifold \Gras ~ consists of the $n$-planes passing
through the origin of $\gc^{n+m}$.
The Pl\"ucker embedding $\iota :$ \Gras$\hookrightarrow {\gcp}^{N(n)-1}$
is given by $\iota (Z)=[Z^{{i_1\ldots i_n}}]\label{embpl}$, where
\begin{equation}
Z=\mb{z}_1\wedge\ldots\wedge\mb{z}_n=
\sum_{1\leq i_1<\ldots <i_n\leq N} Z^{i_1\ldots i_n}
\mb{e}_{i_1}\wedge\ldots\wedge\mb{e}_{i_n},~
\mb{z}_i=\sum_{a=1}^N \hat{Z}_{ia}\mb{e}_a~.   \label{zzz}
\end{equation}
$ Z^{i_1\ldots i_n}$
are the Pl\"ucker coordinates, $\hat{Z}$ =
$(\hat{Z}_{ia})_{1\le i\le n;1\le a\le N}$, $N(n)={N\choose n}=
\displaystyle{\frac{N!}{n!m!}}$ and $\{\mb{e}_1,
...,\mb{e}_N\}$ is a basis of $\gc^N$, $N=n+m$.
Let the vectors  $\mb{z^\sigma_i}$ be such that
$\hat{Z}^{\sigma} \in {\cal V}_\sigma $, where $\sigma$ is a Schubert
symbol. Then
$Z_{\sigma (i)\sigma (\alpha )}~,~i=1,\ldots ,n~ ,
\alpha = n+1,..., N$, are the Pontrjagin coordinates
\begin{equation}
\mb{z^\sigma _i}=\mb{e}_{\sigma (i)}+\sum_{\alpha =n+1}^N
Z_{\sigma (i)\sigma (\alpha )}\mb{e}_{\sigma (\alpha )}~,~i=1,\ldots ,n~ .
\label{16}
\end{equation}

\begin{rem}[Wong\cite{wong}] {\it The cut locus of the point }
{\bf O}$\in$\Gras~
{\it is given by} \label{wogn}
\begin{eqnarray}
\nonumber {\bf CL}_0 & = & \Sigma_0  =	V^m_1 = Z(\omega^m_1)= Z(m-1,m,\ldots ,
m)\\ & =  & \left\{ X\in G_n({\gc}^{n+m})
\vert \dim (X\cap {\bf O}^{\perp})\geq 1\right\}.
\end{eqnarray}
\begin{equation}
V^m_1=\label{74.1}
\cases {{\gcp}^{m-1},& for $n=1 ,$\cr
	      W^m_1\cup W^m_2\cup \ldots W^m_{r-1}\cup W^m_r,&$1<n ,$\cr}
\end{equation}
\begin{equation}
W^m_r=\cases {G_r({\gc}^{\max (m,n)}),&$n\not= m,$\cr
	     {\bf O}^{\perp},&$n=m .$\cr}
\end{equation}
\end{rem}
The following notation is used:
$$V^p_l=\left \{Z\in  G_n({\gc}^{n+m})\vert \dim (Z\cap {\gc}^p)
\geq l\right\}$$
$$\omega^p_l=(\underbrace{p-l,\ldots ,p-l}_l, \underbrace{m,\ldots ,m}_{n-l})$$
$$W^p_l=V^p_l-V^p_{l+1};~~V^p_l=Z(\omega ^p_l);~~ W^p_l=Z'(\omega ^p_l)~$$
$$\omega=\{ 0\leq\omega (1)\leq \ldots \leq\omega (n)\leq m\};~
\sigma (i)=\omega (i)+i,~ i=1,\ldots ,n$$
$$Z(\omega )=\left\{ X\in G_n({\gc}^{n+m}) \vert
 \dim (X\cap {\gc}^{\sigma (i)})\geq i\right\}$$
$$I_{\omega }=\left \{0=i_0<i_1< \ldots <i_{l-1}<i_{l}=n\right\}~$$
$$\omega (i_h)<\omega (i_{h+1}), \omega (i)=\omega (i_{h-1}), i_{h-1}<i
\leq i_h, h=1,\ldots ,l $$
$$Z'(\omega )=\left\{ X\in  G_n({\gc}^{n+m}) \vert
\dim (X\cap {\gc}^{\sigma (i_h)})= i_h,~  i_h\in I_{\omega}\right\} .$$

\begin{rem}{The cut locus for flag manifolds
$G^{\gcm} /P$ has a stratified structure consisting of $r$ $P$-orbits}.
(r=rank).
\end{rem}

\section{\hspace{-4mm}.\hspace{2mm}RESULTS}

\begin{pr}
{If \men~ is an homogeneous algebraic manifold embedded	in a projective
Hilbert space (\ref{emb})
then the polar divisor $\Sigma_0$
can be expressed as $\Sigma_0=\iota^*H_1$, and $\Sigma _0$ is a divisor.}
\end{pr}
{\it Proof.} Use is made of the Cauchy formula	\cite{ber1,ber8}
\begin{equation}\label{cauchy}
(\underline{{\mb e}}_{Z'}, \underline{{\mb e}}_{Z})_{\men}=
( \underline{{\mb e}}_{\iota (Z')},\underline{{\mb e}}_{\iota (Z)})_{\gph},
\end{equation}
where $\iota (Z)=[\underline{{\mb e}}_{Z}]$. We equate with $0$ both sides
of eq. (\ref{cauchy}). The pull-back $\iota^*(H_1)$ of the divisor
$H_1$ is itself a divisor \cite{hartshorne},
because the mapping $\iota$ is an embedding, i.e.
biholomorphic on his image. \pre

\begin{thm}
{\it Let \men~ be a homogeneous simply connected Hodge manifold
admitting the embedding
(\ref{emb}). Let ${\bf M}=\iota^{*}[1]$ be the
unique, up to equivalence, projectively induced line bundle with a given
admissible connection. Then ${\bf M}=[\Sigma_0]$. Moreover, if the homogeneous
manifold \men~ verifies condition B), then ${\bf M} =[{\bf CL}_0]$.
In particular,
the first relation is true for \Ka ian $C$-spaces, while the second one for
hermitian symmetric spaces\/.}
\end{thm}
{\it Proof}. The main part of the proof is based on the following
theorem of Kodaira and Spencer:  For an algebraic manifold there is an
isomorphism of the group ${\mbox{\rm Cl}}(\men )$
of divisor classes with respect to linear equivalence with the Picard group
${\mbox{\rm Pic}}(\men )$, i.e. for
every complex line bundle ${\bf M}$ over an algebraic manifold \men~ there
exists a divisor $D$ such that $[D]={\bf M}$.  The
next ingredient is the following theorem due to Kostant:
 Let \men~ be a simply connected Hodge manifold.
 Then, up to equivalence,
there exists a unique line bundle with	a given curvature matrix
of the hermitian connection, or,
equivalently, with a given admissible connection (Thm. 2.2.1 in \cite{kost}
p. 135). Farther the theorem \ref{cutu} is used.
The information on \Ka ian $C-$ spaces
is extracted from \cite{wolf,ber8}. \pre

We remember also in this context {\it  Bertini theorem:}
{\it Let ${\bf M}$ be a projectively  induced line
bundle over an algebraic manifold \men . Then there is a non-singular divisor
$D$ of \men~ with ${\bf M}=[D]$}. Another formulation reads as follows:
{\it A general hyperplane section $S$ of a connected
non-singular algebraic manifold \men ~in $\gcp^N$ is itself non-singular and for
$n\ge 2$, connected\/} \cite{hartshorne}.

\begin{com}{Generally, the divisor $\Sigma_0$ is singular because it does not
correspond to a general section in Bertini's theorem.}
\end{com}

{\it Proof.} We illustrate the assertion on the pedagogical example
furnished by the Grassmannian $G_2(\gc^4)$.
The Pl\"ucker embedding  is $G_2(\gc^4)\hookrightarrow\gcp^5$.
The coordinate neighbourhoods  ${\cal V}_1-{\cal V}_6$ are presented
in Table 1,
while Table II presents the patches.
%\newpage

In ${\cal V}_1$ the Pl\"ucker coordinates are
$$(p_{12},p_{13},p_{14},p_{23},p_{24},p_{34})=(1,a_{3},a_{4},-a_{1},
-a_{2},a_{1}a_{4}-a_{2}a_{3}).$$
They verify the constrained: $p_{12}p_{34}-p_{13}p_{24}+p_{14}p_{23}=0$.

Let $p_0=(1,0,0,0,0,0)\in{\cal V}_1$ . We want to calculate $\Sigma_0$.
Firstly note that ${\cal V}_1\cap \Sigma_0=\emptyset$.
Then observe that
$\Sigma_0\cap {\cal V}_2=\{b_{3}=0\}$, i.e. an open subset of codimension 1.
We proceed similarly on other coordinate neighbourhoods ${\cal V}_3-
{\cal V}_5$.
On ${\cal V}_6$
$(p_{12},p_{13},p_{14},p_{23},p_{24},p_{34})=(f_{1}f_{4}-f_{3}f_{2},
-f_{3},f_{1},-f_{4},f_{2},1)$. This implies that
$\Sigma_0\cap {\cal V}_6=\{f_{1}f_{4}
-f_{2}f_{3}=0\}$, i.e. disjoint union
of the open subset of codimension 1 already found in ${\cal V}_2-{\cal V}_5$
and the point $(0,0,0,0,0,1)$.
So, in ${\cal V}_6$: $p_{12}=0; p_{14}p_{23}-p_{13}p_{24}=0$.
This is a cone over a quadric surface whose vertex is the point
$(0,0,0,0,0,1)$. The hyperplane $p_{12}=0$ is the
embedded tangent hyperplane of $G_2(\gc ^4)$
of the line $x_1=x_2=0$ in $\gcp^5$. {\it  A general hyperplane
section of $G_2(\gc ^4 )$ is not of the form $p_{12}=0$, since by Bertini's
theorem it has to be smooth}.

In fact, we have also proved the Remark \ref{wogn}
in the particular case
of $G_2(\gc  ^4)$, i.e. $\Sigma_0={\bf CL}_0=V^2_1=W^2_1\cup [(0,0,0,0,0,1)]$,
{\it where $W^2_1$ is a quasiprojective variety of codimension one, while
the point is the singular set.} See details in \cite{arrondo}.\pre

\begin{table}
\caption[Table I]{The Pontrjagin coordinates in the neighbourhoods
${\cal V}_1-{\cal V}_6$ on $G_2(\gc^4)$. }
\begin{center}
\begin{tabular}{|c||c|c|c|c||}\hline
~ & 1 & 2 & 3 & 4  \\ \hline\hline

${\cal V}_1$ & $\begin{array}{c} 1 \\ 0 \end{array}$ &
$ \begin{array}{c} 0 \\ 1 \end{array}$ &
$ \begin{array}{c} a_1 \\ a_3 \end{array}$ &
$ \begin{array}{c} a_2 \\ a_4 \end{array}$ \\ \hline

${\cal V}_2$ & $\begin{array}{c} 1 \\ 0 \end{array}$ &
$ \begin{array}{c} b_1 \\ b_3 \end{array}$ &
$ \begin{array}{c} 0 \\ 1 \end{array}$ &
$ \begin{array}{c} b_2 \\ b_4 \end{array}$ \\ \hline

${\cal V}_3$ & $\begin{array}{c} 1 \\ 0 \end{array}$ &
$ \begin{array}{c} c_1 \\ c_3 \end{array}$ &
$ \begin{array}{c} c_2 \\ c_4 \end{array}$ &
$ \begin{array}{c} 0 \\ 1 \end{array}$ \\ \hline

${\cal V}_4$ & $\begin{array}{c} d_1 \\ d_3 \end{array}$ &
$ \begin{array}{c} 1 \\  0 \end{array}$ &
$ \begin{array}{c} 0 \\ 1 \end{array}$ &
$ \begin{array}{c} d_2 \\ d_4  \end{array}$ \\ \hline

${\cal V}_5$ & $\begin{array}{c} e_1 \\ e_3 \end{array}$ &
$ \begin{array}{c} 1 \\ 0 \end{array}$ &
$ \begin{array}{c} e_2 \\ e_4 \end{array}$ &
$ \begin{array}{c} 0 \\ 1 \end{array}$ \\ \hline

${\cal V}_6$ & $\begin{array}{c} f_1 \\  f_3 \end{array}$ &
$ \begin{array}{c} f_2 \\ f_4 \end{array}$ &
$ \begin{array}{c} 1 \\ 0 \end{array}$ &
$ \begin{array}{c} 0 \\ 1 \end{array}$ \\ \hline\hline

\end{tabular}
\end{center}
\end{table}

\begin{table}
\caption[Table 2]{The change of coordinates on $G_2(\gc^4)$}
\begin{center}
\begin{tabular}{|c||c|c|c|c|c|c||}\hline
~ & ${\cal V}_1$ & ${\cal V}_2$ &
${\cal V}_3$ &	${\cal V}_4$ &	${\cal V}_5$ &	${\cal V}_6$   \\ \hline\hline
${\cal V}_1$ &
$\times$ & $a_3\not= 0 $ & $a_4\not= 0$ & $a_1\not= 0$ & $a_2\not= 0$ &
$\begin{array}{c} a_1a_4-\\
a_2a_3\not= 0\end{array}$\\ \hline

${\cal V}_2$ &
$ b_3\not= 0$ & $\times$ & $b_1\not= 0 $ & $ b_4\not= 0 $ &
$\begin{array}{c} b_1b_4-\\
b_2b_3\not= 0\end{array}$ & $b_2\not= 0$ \\ \hline

${\cal V}_3$ &
$c_3\not= 0$ &$c_4\not= 0$ &  $\times$ & $\begin{array}{c} c_1c_4-\\
c_2c_3\not= 0\end{array}$ & $c_1\not= 0$ &  $c_2\not= 0$
 \\ \hline

${\cal V}_4$ & $ d_3\not= 0$ &$d_1\not= 0$ & $\begin{array}{c} d_1d_4-\\
d_2d_3\not= 0\end{array}$ &
  $\times$ & $d_4\not= 0$ &  $d_2\not= 0$
\\ \hline

${\cal V}_5$ &
$e_3\not= 0$ &	$\begin{array}{c} e_1e_4-\\
e_2e_3\not= 0\end{array}$ &
$e_1\not= 0$ &	$e_4\not= 0$ &
  $\times$ & $e_2\not= 0$
\\ \hline

${\cal V}_6$ &
$\begin{array}{c} f_1f_4-\\
f_2f_3\not= 0\end{array}$ & $f_3\not= 0$ &
$f_1\not= 0$ &	$f_4\not= 0$ &
  $f_2\not= 0$ & $\times$
\\ \hline\hline

\end{tabular}
\end{center}
\end{table}

\vspace*{4ex}
{\bf Acknowledgments}
The author expresses his thanks to Professors
A. Odzijewicz and A.  Strasburger for the possibility
to attend  the XV-th  Workshop on Geometric Methods in Physics.
Discussions during the Workshop with Professors M. Cahen, M. Flato,
Z. Pasternak-Winiarski and M. Schlichenmaier are kindly acknowledged.
Correspondence with Professors J. Burns and E. Arrondo is also
acknowledged. The author is grateful to Professor M. Schlichenmaier
for suggestions in improving the manuscript.

\end{document}